\documentclass[12pt]{article}

\usepackage{amsmath,amsthm,amssymb,stmaryrd,mathrsfs}

\newtheorem{theorem}{Theorem}
\newtheorem{lemma}{Lemma}

\newtheorem{corollary}{Corollary}

\theoremstyle{definition}

\def\fl#1{\left\lfloor#1\right\rfloor}

\renewcommand{\mod}{\text{ mod }}
\def\mand{\qquad \mbox{and} \qquad}

\def\\{\cr}

\def\N{{\mathbb N}}
\def\R{{\mathbb R}}
\def\Z{{\mathbb Z}}

\def\cB{\mathcal B}
\def\cF{\mathcal F}
\def\cI{\mathcal I}
\def\cV{\mathcal V}

\def\mF{\mathscr{F}}

\def\e{{\rm\bf e\/}}
\def\ep{\epsilon}

\def\({\left(}
\def\){\right)}
\def\[{\left[}
\def\]{\right]}
\def\<{\langle}
\def\>{\rangle}
\def\fl#1{\left\lfloor#1\right\rfloor}

\def\le{\leqslant}
\def\ge{\geqslant}

\def \Bt {\cB_{\alpha,\beta}}

\begin{document}

\title{\sc Sums with multiplicative functions
 over a Beatty sequence\footnote{MSC Numbers: 11E25, 11B83.}} % CHECK THIS

\author{
{\sc Ahmet M.~G\"ulo\u glu\footnote{Corresponding author}} \\
{Department of Mathematics} \\
{University of Missouri} \\
{Columbia, MO 65211 USA} \\
{\tt ahmet@math.missouri.edu} \\
\and
{\sc C.~Wesley Nevans} \\
{Department of Mathematics} \\
{University of Missouri} \\
{Columbia, MO 65211 USA} \\
{\tt nevans@math.missouri.edu}}

\date{}

\maketitle

\begin{abstract}
We study sums with multiplicative functions that take values over
a non-homogenous Beatty sequence. We then apply our result in a
few special cases to obtain asymptotic formulas such as the number
of integers in a Beatty sequence representable as a sum of two
squares up to a given magnitude.
\end{abstract}

\newpage

\section{Introduction}
Let $A \ge 1$ be an arbitrary constant, and let $\mF_A$ be the set
of multiplicative functions such that $|f(p)|\le A$ for all primes
$p$, and
\begin{equation} \label{L2norm}
\sum_{n \le N} |f(n)|^2 \le A^2 N\qquad(N\in\N).
\end{equation}
Exponential sums of the form
\begin{equation} \label{ESum}
S_{\alpha,f}(N)=\sum_{n\le N} f(n)
e(n\alpha)\qquad(\alpha\in\R,~f\in\mF_A),
\end{equation}
where $e(z)=e^{2\pi iz}$ for all $z\in\R$, occur frequently in
analytic number theory.  Montgomery and Vaughan have shown
(see~\cite[Corollary 1]{MV}) that the upper bound
\begin{equation} \label{MV Bound}
S_{\alpha,f}(N) \ll_A \frac N {\log N} + \frac{N(\log
R)^{3/2}}{R^{1/2}}
\end{equation}
holds uniformly for all $f \in \mF_A$ provided that $|\alpha-a/q|
\le q^{-2}$ with some reduced fraction $a/q$ for which $2 \le R
\le q \le N/R$. They also proved that this bound is sharp apart
from the logarithmic factor in $R$. In this paper, we use the
Montgomery-Vaughan result to estimate sums of the form
\begin{equation} \label{BSum}
G_{\alpha,\beta,f}(N)=\sum_{\substack{n \le N\\n \in \Bt}} f(n),
\end{equation}
where $\alpha,\beta \in \R$ with $\alpha>1$, $f\in\cF_A$, and
$\Bt$ is the \emph{non-homogenous Beatty sequence} defined by
$$
\cB_{\alpha,\beta}=\big\{n\in\N~:~n=\fl{\alpha m+\beta}\text{~for
some~}m\in\Z\big\}.
$$
Our results are uniform over the family $\mF_A$ and nontrivial
whenever
$$
\lim_{N\to\infty}\frac{\log N}{N\log\log N}\,\biggl|
 \,\sum_{n\le N}f(n)\biggr|=\infty,
$$
a condition which guarantees that the error term in
Theorem~\ref{mainthm} is smaller than the main term. One can
remove this condition, at the expense of losing uniformity with
respect to $f$, and still obtain Theorem \ref{mainthm} for any
bounded arithmetic function $f$ (not necessarily multiplicative)
for which the exponential sums in~\eqref{ESum} satisfy
$$
S_{\alpha,f}(N)=o\(\,\sum_{n \le N}f(n)\)\qquad(N\to\infty).
$$
The general problem of characterizing functions for which this
relation holds appears to be rather difficult; see~\cite{Ba} for
Bachman's conjecture and his related work on this problem.

We shall also assume that $\alpha$ is irrational and of finite
type $\tau$. For an irrational number $\gamma$, the type of
$\gamma$ is defined by
$$
\tau=\sup \bigl\{t\in\R~:~\liminf\limits_{n\to\infty}
~n^t\,\llbracket\gamma n\rrbracket = 0 \bigr\},
$$
where $\llbracket \cdot \rrbracket$ denotes the distance to the
nearest integer. \emph{Dirichlet's approximation theorem}
implies $\tau\ge 1$ for every irrational number $\gamma$.
According to theorems of Khinchin~\cite{Khin} and of
Roth~\cite{Roth1, Roth2}, $\tau=1$ for \emph{almost
all} real numbers (in the sense of the Lebesgue measure) and
\emph{all} irrational algebraic numbers $\gamma$, respectively;
also see~\cite{Bug,Schm}.

Our main result is the following:

\begin{theorem} \label{mainthm}
Let $\alpha,\beta\in\R$ with $\alpha>1$, and suppose that $\alpha$
is irrational and of finite type. Then, for all $f\in\mF_A$ we
have
$$
G_{\alpha,\beta,f}(N)=\alpha^{-1}\sum_{n\le N}
f(n)+O\(\frac{N\log\log N}{\log N}\),
$$
where the implied constant depends only on $\alpha$ and $A$.
\end{theorem}

The following corollaries are immediate applications of
Theorem~\ref{mainthm}:

\begin{corollary} \label{2sqr}
The number of integers not exceeding $N$ that lie in the Beatty
sequence $\cB_{\alpha,\beta}$ and can be represented as a sum of
two squares is 
$$
\#\{n\le N: n\in\cB_{\alpha,\beta},\;n=\square+\square\}=\frac{C N}{\alpha \sqrt{\log N}} + O\(\frac{N\log\log N}{\log N}\)
$$
%asymptotic to $\alpha^{-1}\,C N (\log N)^{-1/2}$
where
\begin{equation} \label{Landau C}
C = 2^{-1/2} \prod_{p \equiv 3 \mod 4} ( 1- p^{-2})^{-1/2}.
\end{equation}
\end{corollary}

To state the next result, we recall that an integer $n$ is said to
be \emph{$k$-free} if $p^k \nmid n$ for every prime $p$.

\begin{corollary} \label{kfree}
For every $k\ge 2$, the number of $k$-free integers not exceeding
$N$ that lie in the Beatty sequence $\cB_{\alpha,\beta}$ is
$$
\#\{n\le N:n\in\cB_{\alpha,\beta},\;n \text{ is }k\text{-free} \}=\alpha^{-1}\zeta^{-1}(k)N+O\(\frac{N\log\log N}{\log N}\)
$$
%asymptotic to $\alpha^{-1}\,\zeta^{-1}(k)N$ 
where $\zeta(s)$ is the Riemann zeta function.
\end{corollary}

Finally, we consider the average value of the number of
representations of an integer from a Beatty sequence as a sum of
four squares. Our result is the following:

\begin{corollary} \label{4sqr}
Let $r_4(n)$ denote the number of representations of $n$ as a sum
of four squares. Then,
$$
\sum_{\substack{n\le N\\n\in \Bt}}r_4(n)=\frac{\pi^2N^2}{2\alpha}
+O\(\frac{N^2\log\log N}{\log N}\),
$$
where the implied constant depends only on $\alpha$.
\end{corollary}

Any implied constants in the symbols $O$ and $\ll$ may depend on the parameters $\alpha$ and $A$ but are absolute otherwise. We recall that the notation $X\ll Y$ is equivalent to $X=O(Y)$.

\bigskip

\noindent {\bf Acknowledgments.} We would like to thank William
Banks and Igor Shparlinski for their helpful comments and careful reading of the
original manuscript.

\section{Preliminaries}
\subsection{Discrepancy of fractional parts}
\label{sec:discr}

We define the \emph{discrepancy} $D(M)$ of a sequence of 
real numbers $b_1,b_2,\ldots,b_M\in[0,1)$ by
\begin{equation}
\label{eq:descr_defn}
D(M)=\sup_{\cI\subseteq[0,1)}\left|\frac{V(\cI,M)}{M}-|\cI|\,\right|,
\end{equation}
where the supremum is taken over all possible subintervals $\cI=(a,c)$ of the
interval $[0, 1)$, $V(\cI,M)$ is the number of positive integers
$m\le M$ such that $b_m\in\cI$, and $|\cI|=c-a$ is the length of $\cI$.

If an irrational number $\gamma$ is of finite type, we let $D_{\gamma,\delta}(M)$ denote the discrepancy of the sequence of fractional parts 
$(\{\gamma m + \delta\})_{m=1}^M$. By~\cite[Theorem~3.2, Chapter~2]{KuNi} we have:

\begin{lemma}
\label{lem:discr_with_type}  For a fixed irrational number $\gamma$ of finite type $\tau$ and for all $\delta\in\R$ we have:
$$
D_{\gamma,\delta}(M)\le M^{-1/\tau+o(1)}\qquad(M\to\infty),
$$
where the function defined by $o(\cdot)$ depends only on $\gamma$.
\end{lemma}

\subsection{Numbers in a Beatty sequence}

The following is standard in characterizing the elements of the Beatty sequence $\Bt$:

\begin{lemma} \label{lem:Beatty_values}
Let $\alpha,\beta\in\R$ with $\alpha>1$, and set $\gamma=\alpha^{-1}$, $\delta=\alpha^{-1}(1-\beta)$. Then,
$n=\left\lfloor\alpha m+\beta\right\rfloor$ for some $m\in\Z$ if and only if $0<\{\gamma n+\delta\}\le\gamma$.
\end{lemma}
From Lemma~\ref{lem:Beatty_values}, an integer $n$ lies in
$\Bt$ if and only if $n\ge 1$ and $\psi(n)=1$,
where $\psi$ is the periodic function with period one whose values
on the interval $(0,1]$ are given by
$$
\psi(x) = \left\{  \begin{array}{ll}1& \quad \hbox{if $0<x\le \gamma$}; \\0& \quad \mbox{if $\gamma<x\le 1$}.\end{array} \right.
$$
We wish to approximate $\psi$ by a function whose Fourier series representation is well behaved. This will gives rise
to the fore mentioned exponential sum $S_{\alpha,f}(N)$. To this end we use the result of Vinogradov (see~\cite[Chapter~I, Lemma~12]{Vin}), which states that for any $\Delta$ such that
$$ 0 < \Delta < \frac{1}{8} \mand \Delta\le\frac{1}{2}\min\{\gamma,1-\gamma\},
$$
there exists a real-valued function $\Psi$ with the following properties:
\begin{itemize}
\item[$(i)$~~] $\Psi$ is periodic with period one;

\item[$(ii)$~~] $0 \le\Psi(x)\le 1$ for all $x\in\R$;

\item[$(iii)$~~] $\Psi(x)=\psi(x)$ if $\Delta\le \{x\}\le
\gamma-\Delta$ or if $\gamma+\Delta\le \{x\}\le 1-\Delta$;

\item[$(iv)$~~] $\Psi$ can be represented by a Fourier series:
$$
\Psi(x)=\sum_{k\in\Z}g(k)\e(kx),
$$
where $g(0)=\gamma$, and the Fourier coefficients satisfy the
uniform bound
\begin{equation} \label{eq:coeffbounds}
g(k) \ll \min \big\{ |k|^{-1}, |k|^{-2} \Delta^{-1} \big\}
\qquad(k\ne 0).
\end{equation}
\end{itemize}

\section{Proofs}

\subsection{Proof of Theorem \ref{mainthm}}

Using Lemma~\ref{lem:Beatty_values}, we rewrite the
sum~\eqref{BSum} in the form
$$
G_{\alpha,\beta,f}(N)= \sum_{n\le N} f(n) \psi(\gamma n+\delta).
$$
Replacing $\psi$ by $\Psi$ we have
\begin{equation}
\label{eq:one} G_{\alpha,\beta,f}(N)= \sum_{n \le N} f(n) \Psi(
\gamma n+\delta) + O\(\sum_{n\in \cV(\cI,N)}f(n)\),
\end{equation}
where $\cV(\cI,N)$ is the set of positive integers $n\le N$ for
which
$$ \{\gamma n+\delta\} \in \cI=
[0,\Delta)\cup(\gamma-\Delta,\gamma+\Delta) \cup(1-\Delta,1). $$
Since $|\cI|=4\Delta$, it follows from
Lemma~\ref{lem:discr_with_type} and the
definition~\eqref{eq:descr_defn} that
$$
\big|\cV(\cI,N)\big|\ll\Delta N+N^{1-1/(2\tau)},
$$
where we have used the fact that $\alpha$ and $\gamma$ have the
same type $\tau$. Thus, taking~\eqref{L2norm} into account, we
have by the Cauchy inequality:
\begin{equation}
\label{eq:two}
\begin{split}
\Biggl|\,
\sum_{n\in\cV(\cI,N)}f(n)\Biggr|&\le\big|\cV(\cI,N)\big|^{1/2}
\(\,\sum_{n\le N}|f(n)|^2\)^{1/2}\\
&\ll\Delta^{1/2}N+N^{1-1/(4\tau)}.
\end{split}
\end{equation}

Next, let $K\ge\Delta^{-1}$ be a large real number (to be
specified later), and let $\Psi_K$ be the trigonometric polynomial
given by
\begin{equation}
\label{eq:PKdefn} \Psi_K(x)=\sum_{|k|\le
K}g(k)\e(kx)=\gamma+\sum_{0<|k|\le K}g(k)\e(kx)\qquad(x\in\R).
\end{equation}
Using \eqref{eq:coeffbounds} we see that the estimate
$$
\Psi_K(x)=\Psi(x)+O(K^{-1}\Delta^{-1})
$$
holds uniformly for all $x\in\R$; therefore,
\begin{equation}
\label{eq:three} \sum_{n\le N}f(n)\psi(\gamma n+\delta)=\sum_{n\le
N}f(n)\Psi_K(\gamma n+\delta)+O\(K^{-1}\Delta^{-1}N\),
\end{equation}
where we have used the bound $\sum_{n\le N}|f(n)|\ll N$, which
follows from~\eqref{L2norm}.

Combining~\eqref{eq:one}, \eqref{eq:two}, \eqref{eq:PKdefn}
and~\eqref{eq:three} we derive that
$$
G_{\alpha,\beta,f}(N)=\gamma \sum_{n \le N} f(n) + H(N)+
O\(K^{-1}\Delta^{-1}N +\Delta^{1/2}N+N^{1-1/(4\tau)}\),
$$
where
$$
H(N) = \sum_{0 < |k| \le K} g(k) \e(k\delta) S_{k\gamma,f}(N).
$$
Put $R=(\log N)^3$. We claim that, if $N$ is sufficiently large,
then for every $k$ in the above sum there is a reduced fraction
$a/q$ such that $|\alpha-a/q| \le q^{-2}$ and $R \le q \le N/R$.
Assuming this for the moment, \eqref{MV Bound} implies that
$$
S_{k\gamma,f}(N)\ll\frac{N}{\log N}\qquad(0<|k|\le K),
$$
and using~\eqref{eq:coeffbounds} we deduce that
$$
H(N)\ll\frac{N \log K}{\log N}.
$$
Therefore,
$$
G_{\alpha,\beta,f}(N)-\gamma \sum_{n \le N} f(n)\ll \frac{N \log
K}{\log N}+K^{-1}\Delta^{-1}N +\Delta^{1/2}N.
$$
To balance the error terms, we choose
$$
\Delta=(\log N)^{-2}\mand K=\Delta^{-3/2}=(\log N)^3,
$$
obtaining the bound stated in the theorem.

To prove the claim, let $k$ be an integer with $0<|k|\le K=(\log
N)^3$, and let $r_i=a_i/q_i$ be the $i$-th convergent in the
continued fraction expansion of~$k\gamma$. Since $\gamma$ is of
finite type $\tau$, for every $\ep>0$ there is a constant
$C=C(\gamma,\ep)$ such that
$$
C(|k|q_{i-1})^{-(\tau+\ep)}<\llbracket \gamma |k| q_{i-1}
\rrbracket \le \bigl| \gamma |k|q_{i-1}-a_{i-1} \bigr| \le
q_i^{-1}.
$$
Put $\ep=\tau$, and let $j$ be the least positive integer for
which $q_j\ge R$ (note that $j\ge 2$). Then,
$$
R\le q_j\ll (|k|q_{i-1})^{2\tau}\le (KR)^{2\tau}=(\log N)^{6\tau},
$$
and it follows that $R\le q_j\le N/R$ if $N$ is sufficiently
large, depending only on $\alpha$. This concludes the proof.

\subsection{Proof of Corollary~\ref{2sqr}}

Let $f(n)$ be the characteristic function of the set of integers
that can be represented as a sum of two squares.  Then
Corollary~\ref{2sqr} follows immediately from
Theorem~\ref{mainthm} and the
asymptotic formula (see for example~\cite{Shan}): %CWN-replaced the line so that we have an asymptotic with error term
%well known result ofLandau~\cite{La} which asserts that
$$
\sum_{n\le N} f(n)=\frac{CN}{(\log N)^{1/2}} + O\( \frac N {(\log N)^{3/2}}\),
$$
where $C$ is given by \eqref{Landau C}.

\subsection{Proof of Corollary~\ref{kfree}}

Fix $k\ge 2$ and let $f(n)$  be the characteristic function of the
set of $k$-free integers.  Then Corollary~\ref{kfree} follows from
Theorem~\ref{mainthm} and the following estimate of
Gegenbauer~\cite{Ge} for the number of $k$-free integers not
exceeding $N$:
$$
\sum_{n\le N}f(N)=\zeta^{-1}(k)N+O\(N^{1/k}\).
$$

\subsection{Proof of Corollary~\ref{4sqr}}

Put $f(n)=r_4(n)/(8n)$. {From} Jacobi's formula for $r_4(n)$,
namely
$$
r_4(n)=8(2+(-1)^n) \sum_{\substack{d\,\mid\,n\\d\text{~odd}}}d
\qquad(n\ge 1),
$$
it follows that $f(n)$ is multiplicative, and $f(p)\le 3/2$ for
every prime~$p$. Moreover, using the formula of
Ramanujan~\cite{Ram} (see also~\cite{Smith}):
$$
\sum_{n\le N}\sigma^2(n)=\frac56\,\zeta(3)N^3+O(N^2(\log N)^2),
$$
we have by partial summation:
\begin{eqnarray*}
\sum_{n \le N} |f(n)|^2 \le \sum_{n \le N} \frac{\sigma^2
(n)}{n^2} = \frac52\,\zeta(3)N+ O((\log N)^3).
\end{eqnarray*}
Therefore, $f(n)\in\mF_A$ for some constant $A\ge 1$. Applying
Theorem~\ref{mainthm}, we deduce that
$$
\sum_{\substack{n \le N\\n \in
\Bt}}\frac{r_4(n)}{n}=\alpha^{-1}\sum_{n\le N}
\frac{r_4(n)}{n}+O\(\frac{N\log\log N}{\log N}\),
$$
where the implied constant depends only on $\alpha$.

{From} the asymptotic formula (see for example~\cite[p22]{IwKow}):
$$
\sum_{n\le N} r_4 (n) = \frac{\pi^2N^2}{2}+O(N\log N),
$$
we have by partial summation:
\begin{eqnarray*}
\sum_{n\le N}\frac{r_4(n)}{n}=\pi^2N+O((\log N)^2).
\end{eqnarray*}
Consequently,
$$
\sum_{\substack{n\le
N\\n\in\Bt}}\frac{r_4(n)}{n}=\alpha^{-1}\pi^2N+O\(\frac{N\log\log
N}{\log N}\).
$$
Using partial summation once more, we obtain the statement of
Corollary~\ref{4sqr}.

\end{document}